\input amstex
\input Amstex-document.sty
\loadbold
\font\teneuf=eufm10 \font\seveneuf=eufm7 \font\fiveeuf=eufm5
\font\tenmsy=msbm10 \font\sevenmsy=msbm7 \font\fivemsy=msbm5
\font\tenmsx=msam10 \font\sixmsx=msam6 \font\fivemsx=msam5
\textfont7=\teneuf \scriptfont7=\seveneuf
\scriptscriptfont7=\fiveeuf
\textfont8=\tenmsy \scriptfont8=\sevenmsy
\scriptscriptfont8=\fivemsy
\textfont9=\tenmsx \scriptfont9=\sixmsx
\scriptscriptfont9=\fivemsx
\def\gr{\fam7 \teneuf}

\def\grm{{\gr m}} \def\grM{{\gr M}} \def\grS{{\gr S}}
\def\grB{{\gr B}}

\pageno 207

\def\shorttitle{Exponential Sums and Diophantine Problems}
\def\shortauthor{T. D. Wooley}

\topmatter %
\title\nofrills{\boldHuge Diophantine Methods for Exponential Sums, and
Exponential Sums for Diophantine Problems}
\endtitle

\author \Large Trevor D. Wooley* \endauthor

\thanks *Department of Mathematics, University of Michigan, East Hall,
525 East University Avenue, Ann Arbor, MI 48109-1109, USA. E-mail:
wooley\@umich.edu \endthanks

\abstract\nofrills \centerline{\boldnormal Abstract}

\vskip 4.5mm

{\ninepoint Recent developments in the theory and application of the Hardy- \linebreak Littlewood method are
discussed, concentrating on aspects associated with diagonal diophantine problems. Recent efficient differencing
methods for estimating mean values of exponential sums are described first, concentrating on developments
involving smooth Weyl sums. Next, arithmetic variants of classical inequalities of Bessel and Cauchy-Schwarz are
discussed. Finally, some emerging connections between the circle method and arithmetic geometry are mentioned.

\vskip 4.5mm


}
\endabstract
\endtopmatter

\document

\baselineskip 4.5mm \parindent 8mm

\vskip 8mm

\specialhead \noindent \boldLARGE 1. Introduction \endspecialhead

Over the past fifteen years or so, the Hardy-Littlewood method has experienced a renaissance that has left
virtually no facet untouched in its application to diophantine problems. Our purpose in this paper is to sketch
what might be termed the past, present, and future of these developments, concentrating on aspects associated with
diagonal diophantine problems, and stressing modern developments that make increasing use of less traditional
diophantine input within ambient analytic methods. We avoid discussion of the Kloosterman method and its important
recent variants (see [5] and [8]), because
 the underlying ideas seem inherently constrained to quadratic, and
occasionally cubic, diophantine problems. Our account begins with a brief
introduction to the Hardy-Littlewood (circle) method, using Waring's problem
as the basic example. The discussion here illustrates well the issues involved
in the analysis of systems of diagonal equations over arbitrary algebraic
extensions of ${\Bbb Q}$, and motivates that associated with more general
systems of homogeneous equations (see [1] and [14]).

Let $s$ and $k$ be natural numbers with $s>k\ge 2$, and consider an integer $n$
 sufficiently large in terms of $s$ and $k$. The circle method employs Fourier
analysis in order to obtain asymptotic information concerning the number,
$R(n)=R_{s,k}(n)$, of integral solutions of the equation $x_1^k+\dots
+x_s^k=n$. Write $P=n^{1/k}$ and define the exponential sum
$f(\alpha )=f(\alpha ;P)$ by
$$f(\alpha )=\sum_{1\le x\le P}e(\alpha x^k),$$
wherein $e(z)$ denotes $e^{2\pi iz}$. Then it follows from orthogonality that
$$R(n)=\int_0^1f(\alpha )^se(-n\alpha )d\alpha .$$
When $\alpha $ is well-approximated by rational numbers with small
denominators, one has sharp asymptotic information concerning $f(\alpha )$. In
order to be precise, let $Q$ satisfy $1\le Q\le {1\over 2}P^{k/2}$, and define
the {\it major arcs} $\grM =\grM (Q)$ to be the union of the intervals
$\grM (q,a)=\{ \alpha \in [0,1): |q\alpha -a|\le QP^{-k}\} $, with $0\le a\le
q\le Q$ and $(a,q)=1$. Also, put
$$S(q,a)=\sum_{r=1}^qe(ar^k/q)\quad \text{and}\quad v(\beta )=\int_0^P
e(\beta \gamma^k)d\gamma ,$$
and define $f^*(\alpha )$ for $\alpha \in [0,1)$ by taking
$f^*(\alpha )=q^{-1}S(q,a)v(\alpha -a/q)$, when $\alpha $ lies in $\grM (q,a)
\subseteq \grM (Q)$, and otherwise by setting $f^*(\alpha )=0$. Then the
sharpest available estimate (see Theorem 4.1 of [16]) establishes that
\footnote{Given a complex-valued function $f(t)$ and positive function $g(t)$,
we use Vinogradov's notation $f(t)\ll g(t)$, or Landau's notation
$f(t)=O(g(t))$, to mean that when $t$ is large, there is a positive number $C$
for which $f(t)\le Cg(t)$. Similarly, we write $f(t)\gg g(t)$ when $g(t)\ll
f(t)$, and $f(t)\asymp g(t)$ when $f(t)\ll g(t)\ll f(t)$. Also, we write
$f(t)=o(g(t))$ when as $t\rightarrow \infty $, one has $f(t)/g(t)\rightarrow
0$. Finally, we use the convention that whenever $\epsilon $ occurs in a
formula, then it is asserted that the statement holds for each fixed positive
number $\epsilon $.}
$f(\alpha )=f^*(\alpha )+O(Q^{1/2+\epsilon })$, uniformly for $\alpha \in \grM
(Q)$. The functions $S(q,a)$ and $v(\beta )$ are rather well-understood, and
thus one deduces that whenever $s\ge \max \{ 4,k+1\} $ and $Q\le P$, then
$$\int_\grM f(\alpha )^se(-n\alpha )d\alpha ={{\Gamma (1+1/k)^s}\over {\Gamma
(s/k)}}\grS_{s,k}(n)n^{s/k-1}+O(n^{s/k-1-\delta }),\tag1.1$$
for a suitable positive number $\delta $. Here, the $\Gamma $-function is that
familiar from classical analysis, and the
{\it singular series} $\grS_{s,k}(n)$ is equal to the product of $p$-adic
densities $\prod_pv_p(n)$, where for each prime $p$ we write
$$v_p(n)=\lim_{h\rightarrow \infty }p^{h(1-s)}\text{card} \{ {\bold x}\in
({\Bbb Z}/p^h{\Bbb Z})^s\, :\, x_1^k+\dots +x_s^k\equiv n\pmod{p^h}\} .$$

An asymptotic formula for $R(n)$, with leading term determined by the major arc
 contribution (1.1), now follows provided that the corresponding contribution
arising from the minor arcs $\grm =[0,1)\setminus \grM $ is asymptotically
smaller. Although such is conjectured to hold as soon as $s\ge \max
\{ 4,k+1\} $, this is currently known only for larger values of $s$. It is here
 that energy is focused in current research. One typically estimates the minor
arc contribution via an inequality of the type
$$\left| \int_\grm f(\alpha )^se(-n\alpha )d\alpha \right| \le \left(
\sup_{\alpha \in \grm }|f(\alpha )|\right)^{s-2t}\int_0^1|f(\alpha )|^{2t}
d\alpha .\tag1.2$$
For suitable choices of $t$ and $Q$, one now seeks bounds of the shape
$$\sup_{\alpha \in \grm }|f(\alpha )|\ll P^{1-\tau +\epsilon }\quad \text{and}
\quad \int_0^1|f(\alpha )|^{2t}d\alpha \ll P^{2t-k+\delta +\epsilon },\tag1.3$$
with $\tau >0$ and $\delta $ small enough that $(s-2t)\tau >\delta $. The right
 hand side of (1.2) is then $o(n^{s/k-1})$, which is smaller than the main term
 of (1.1) whenever $\grS_{s,k}(n)\gg 1$. The latter is assured provided that
non-singular $p$-adic solutions can be found for each prime $p$, and in any
case when $s\ge 4k$. Classically, one has two apparently
incompatible approaches toward establishing the estimates (1.3). On one side
is the differencing approach introduced by Weyl [23], and pursued by Hua [9],
that yields an asymptotic formula for $R(n)$ whenever $s\ge 2^k+1$. The ideas
introduced by Vinogradov [21], meanwhile, provide the desired asymptotic
formula when $s>Ck^2\log k$, for a suitable positive constant $C$.

\vskip 8mm

\specialhead \noindent \boldLARGE 2. Efficient differencing and smooth Weyl sums\endspecialhead

Since the seminal work of Vaughan [15], progress on diagonal diophantine problems has been based, almost
exclusively, on the use of smooth numbers, by which we mean integers free of large prime factors. In brief, one
seeks serviceable substitutes for the estimates (1.3) with the underlying summands restricted to be smooth, the
hope being that this restriction might lead to sharper bounds. Before describing the kind of conclusions now
available, we must introduce some notation. Let ${\Cal A}(P,R)$ denote the set of natural numbers not exceeding
$P$, all of whose prime divisors are at most $R$, and define the associated exponential sum $h(\alpha )=h(\alpha
;P,R)$ by
$$h(\alpha ;P,R)=\sum_{x\in {\Cal A}(P,R)}e(\alpha x^k).$$
When $t$ is a positive integer, we consider the mean value
$S_t(P,R)=\int_0^1|h(\alpha )|^{2t}d\alpha $, which, by orthogonality, is
equal to the number of solutions of the diophantine
 equation $x_1^k+\dots +x_t^k=y_1^k+\dots +y_t^k$, with $x_i,y_i\in
{\Cal A}(P,R)$ $(1\le i\le t)$. We take $R\asymp P^\eta $ in the ensuing
discussion, with $\eta $ a small positive number\footnote{We adopt the
convention that whenever $\eta $ appears in a statement, implicitly or
explicitly, then it is asserted that the statement holds whenever $\eta >0$ is
sufficiently small in terms of $\epsilon $.}. In these circumstances one has
$\text{card}({\Cal A}(P,R))\sim c(\eta )P$, where the positive number
$c(\eta )$ is given by the Dickman function, and it follows that $S_t(P,R)\gg
P^t+P^{2t-k}$. It is conjectured that in fact $S_t(P,R)\ll P^\epsilon
(P^t+P^{2t-k})$. We refer to the exponent $\lambda_t$ as being {\it
permissible} when, for each $\epsilon >0$, there exists a positive number
$\eta =\eta (t,k,\epsilon )$ with the property that whenever $R\le P^\eta $,
one has $S_t(P,R)\ll P^{\lambda_t+\epsilon }$. One expects that the exponent
$\lambda_t=\max \{ t,2t-k\} $ should be permissible, and with this in mind we
say that $\delta_t$ is an {\it associated exponent} when
$\lambda_t=t+\delta_t$ is permissible, and that $\Delta_t$ is an
{\it admissible exponent} when $\lambda_t=2t-k+\Delta_t$ is permissible.

The computations required to determine sharp permissible exponents for a
specific value of $k$ are substantial (see [20]), but for larger $k$ one may
summarise some general features of these exponents. First, for $0\le t\le 2$
and $k\ge 2$, it is essentially classical that the exponent $\delta_t=0$ is
associated, and recent work of Heath-Brown [6] provides the same conclusion
also when $t=3$ and $k\ge 238,607,918$. When $t=o(\sqrt{k})$, one finds that
associated exponents exhibit {\it quasi-diagonal behaviour}, and satisfy the
property that $\delta_t\rightarrow 0$ as $k\rightarrow \infty $. To be precise,
 Theorem 1.3 of [28] shows that whenever $k\ge 3$ and $2<t\le 2e^{-1}k^{1/2}$,
then the exponent
$$\delta_t={{4k^{1/2}}\over {et}}\exp \left( -{{4k}\over {e^2t^2}}\right) ,
\tag2.1$$
is associated. For larger $t$, methods based on repeated efficient differencing
 yield the sharpest estimates. Thus, the corollary to Theorem 2.1 of [26]
establishes that for $k\ge 4$, an admissible exponent $\Delta_t$ is given by
 the positive solution of the equation $\Delta_te^{\Delta_t/k}=ke^{1-2t/k}$.
The exponent $\lambda_t=2t-k+ke^{1-2t/k}$ is therefore always permissible.
Previous to repeated efficient differencing, analogues of these permissible
exponents had a term of size $ke^{-t/k}$ in place of $ke^{1-2t/k}$ (see [15]),
so that in a sense, the modern theory is twice as powerful as that available
hitherto.

The above discussion provides a useable analogue of the mean-value estimate in
(1.3). We turn next to localised minor arc estimates. Take $Q=P$, and define
$\grm $ as in the introduction. Suppose that $s$, $t$ and $w$ are parameters
with $2s\ge k+1$ for which $\Delta_s$, $\Delta_t$ and $\Delta_w$ are admissible
 exponents, and define
$$\sigma (k)={{k-\Delta_t-\Delta_s\Delta_w}\over {2(s(k+\Delta_w-\Delta_t)+
tw(1+\Delta_s))}}.$$
Then Corollary 1 to Theorem 4.2 of [27] shows that $\sup_{\alpha \in \grm
}|h(\alpha )|\ll P^{1-\sigma (k)+\epsilon }$, and for large $k$ this estimate
holds with $\sigma (k)^{-1}=k(\log k+O(\log \log k))$. Applying an analogue of
(1.2) with $h$ in place of $f$, and taking
\footnote{We write $[z]$ to denote $\max \{ n\in {\Bbb Z}\, :\, n\le z\} $.}
$t=[{1\over 2}k(\log k+\log \log k+1)]$ and $s=2t+k+[Ak\log \log k/\log k]$,
for a suitable $A>0$, we deduce from our discussion of permissible exponents
that $\int_\grm h(\alpha )^se(-n\alpha )d\alpha =o(n^{s/k-1})$. By considering
the representations of a given integer $n$ with all of the $k$th powers
$R$-smooth, it is now apparent that a modification of the argument sketched in
the introduction shows that $R(n)\gg \grS_{s,k}(n)n^{s/k-1}$ as soon as one
confirms that
$$\int_\grM h(\alpha )^se(-n\alpha )d\alpha \sim c(\eta )^s{{\Gamma (1+1/k)^s}
\over {\Gamma (s/k)}}\grS_{s,k}(n)n^{s/k-1}.\tag2.2$$
Sharp asymptotic information concerning $h(\alpha )$ is available throughout
$\grM (Q)$ only when $Q$ is a small power of $\log P$, and so the proof of
(2.2) involves pruning technology. Such machinery, in this case designed to
estimate the contribution from a set of the shape $\grM (P)\setminus \grM
((\log P)^\delta )$, has evolved into a powerful tool. Such issues can be
handled these days with a number of variables barely exceeding
$\max \{ 4,k+1\} $.

This approach leads to the best known upper bounds on the function $G(k)$ in
Waring's problem, defined to be the least integer $r$ for which all
sufficiently large natural numbers are the sum of at most $r$ positive integral
 $k$th powers.

{\bf Theorem 2.1. } \it One has $G(k)\le k(\log k+\log \log k+2+O(\log \log
k/\log k))$.
\rm

This upper bound (Theorem 1.4 of [27]) refines an earlier one of asymptotically
 similar strength (Corollary 1.2.1 of [24]) that gave the first sizeable
improvement of Vinogradov's celebrated bound $G(k)\le (2+o(1))k\log k$, dating
from 1959 (see [22]). Aside from Linnik's bound $G(3)\le 7$ (see [11]), all of
the sharpest known bounds on $G(k)$ for smaller $k$ are established using
variants of these methods. Thus one has $G^\# (4)\le 12$ (see [15], and here
the $\# $ denotes that there are congruence conditions modulo $16$), $G(5)\le
17$, $G(6)\le 24$, $G(7)\le 33$, $G(8)\le 42$, $G(9)\le 50$, $G(10)\le 59$,
$G(11)\le 67$, $G(12)\le 76$, $G(13)\le 84$, $G(14)\le 92$, $G(15)\le 100$,
$G(16)\le 109$, $G(17)\le 117$, $G(18)\le 125$, $G(19)\le 134$, $G(20)\le 142$
(see [17], [18], [19], [20]).

Unfortunately, shortage of space obstructs any but the crudest account of the
ideas underlying the proof of the mean value estimates that supply the above
permissible exponents. The use of exponential sums over smooth numbers occurs
already in work of Linnik and Karatsuba (see [10]), but only with Vaughan's new
 iterative method [15] is a flexible homogeneous approach established. An
alternative formulation suitable for repeated efficient differencing is
introduced by the author in [24]. Suppose that the exponent $\lambda_s$ is
permissible, and consider a polynomial $\psi \in {\Bbb Z}[t]$ of degree
$d\ge 2$. Given positive numbers $M$ and $T$ with $M\le T$, and an element
$x\in {\Cal A}(T,R)$ with $x>M$, there exists an integer $m$ with $m\in
[M,MR]$ for which $m|x$. Consequently, by applying a {\it fundamental lemma} of
 combinatorial flavour, one may bound the number of integral solutions of the
equation
$$\psi (z)-\psi (w)=\sum_{i=1}^s(x_i^k-y_i^k),\tag2.3$$
with $1\le z,w\le P$ and $x_i,y_i\in {\Cal A}(T,R)$ $(1\le i\le s)$, in terms
of the number of integral solutions of the equation
$$\psi (z)-\psi (w)=m^k\sum_{i=1}^s(u_i^k-v_i^k),\tag2.4$$
with $1\le z,w\le P$, $M<m\le MR$, $(\psi'(z)\psi'(w),m)=1$ and $u_i,v_i\in
{\Cal A}(T/M,R)$ $(1\le i\le s)$. The implicit congruence condition $\psi
(z)\equiv \psi (w)\pmod{m^k}$ may be analytically refined to the stronger one
$z\equiv w\pmod{m^k}$, and in this way one is led to replace the expression
$\psi (z)-\psi (w)$ by the difference polynomial
$\psi_1(z;h;m)=m^{-k}(\psi (z+hm^k)-\psi (z))$. Notice that when $M\ge
P^{1/k}$, one is forced to conclude that $z=w$, and then the number of
solutions of (2.4) is bounded above by $PMRS_s(T/M,R)\ll P^{1+\epsilon
}M(T/M)^{\lambda_s}$. Otherwise, following an application of Schwarz's
inequality to the associated mean value of exponential sums, one may recover an
 equation of the shape (2.3) in which $\psi (z)$ is replaced by $\psi_1(z)$,
and $T$ is replaced by $T/M$, and repeat the process once again. This gives a
repeated differencing process that hybridises that of Weyl with the ideas of
Vinogradov.

It is now possible to describe a strategy for bounding a permissible exponent
$\lambda_{s+1}$ in terms of a known permissible exponent $\lambda_s$. We
initially take $T=P$ and $\psi (z)=z^k$, and observe that $S_{s+1}(P,R)$ is
bounded above by the number of solutions of (2.3). We apply the above efficient
 differencing process successively with appropriate choices for $M$ at each
stage, say $M=P^{\phi_i}$, with $0\le \phi_i\le 1/k$, for the $i$th
differencing operation. After some number of steps, say $j$, we take
$\phi_j=1/k$ in order to force the above diagonal situation that is easily
estimated. One then optimises choices for the $\phi_i$ in order to extract the
sharpest upper bound for $S_{s+1}(P,R)$, and this in turn yields a permissible
exponent $\lambda_{s+1}$. It transpires that in this simplified treatment,
successive admissible exponents are related by the formula
$\Delta_{s+1}=\Delta_s(1-\phi )+k\phi -1$, wherein one may take $\phi $ very
close to $1/(k+\Delta_s)$. Thus one finds that $\Delta_{s+1}$ is essentially
$\Delta_s(1-2/(k+\Delta_s))$, an observation that goes some way to explaining
how it is that this method is about twice as strong as previous approaches
that would correspond to choices of $\phi $ close to $1/k$.

Refined versions of this differencing process make use of all known permissible
 exponents $\lambda_s$ in order to estimate a particular exponent $\lambda_t$,
and in such circumstances the process becomes highly iterative, and entails
significant computation. Such variants make use of refined Weyl estimates for
difference polynomials, and estimates for the number of integral points on
curves and surfaces (see [20]). Variants of these methods apply also in the
situation of Vinogradov's mean value theorem (see [25]), smooth Weyl sums with
polynomial arguments (see [29]), and even for sums relevant to counting
rational lines on hypersurfaces (see [12]).

Frequent reference to underlying diophantine equations seems to limit these
methods to estimating even moments of smooth Weyl sums, and until recently
fractional moments could be estimated only by applying H\"older's inequality to
 interpolate linearly between permissible exponents. However, a method [28] is
now available that permits fractional moments to be estimated non-trivially,
thereby ``breaking classical convexity'', and moreover the number of variables
being differenced need not even be an integer. These new estimates can be
applied to sharpen permissible exponents (with integral argument), and indeed
the associated exponent (2.1) is established in this way. Another consequence
[32] of these developments is the best available lower bound for $N(X)$, which
we define to be the number of integers not exceeding $X$ that are represented
as the sum of three positive integral cubes. One has $N(X)\gg
X^{1-\xi/3-\epsilon }$, where $\xi =(\sqrt{2833}-43)/41=0.24941301\dots $
arises from the permissible exponent $\lambda_3=3+\xi$ for $k=3$. Earlier,
Vaughan [15] obtained an estimate of the latter type with $13/4$ in place of
$3+\xi $.

\vskip 8mm

\specialhead \noindent \boldLARGE 3. Arithmetic variants of Bessel's inequality\endspecialhead

Already in our opening paragraph we alluded to some of the applications accessible to the methods of \S2. We now
turn to less obvious applications that have experienced recent progress. We illustrate ideas once again with a
simple example, and consider the set ${\Cal Z}(N)$ of integers $n$, with $N/2<n\le N$, that are {\it not}
represented as the sum of $s$ positive integral $k$th powers. The standard approach to estimating
$Z(N)=\text{card}({\Cal Z}(N))$ is via Bessel's inequality. We now take $P=N^{1/k}$. When $\grB \subseteq [0,1)$,
write $R^*(n;\grB )=\int_\grB h( \alpha )^se(-n\alpha )d\alpha $, and write also $R^*(n)=R^*(n;[0,1))$. The theory
of \S2 ensures that when $Q$ is a sufficiently small power of $\log P$, and $s\ge 4k$, then $R^*(n;\grM )\asymp
n^{s/k-1}$. Under such circumstances, an application of Bessel's inequality reveals that $Z(N)$ is bounded above
by
$$\align \sum_{N/2<n\le N}\left| {{R^*(n)-R^*(n;\grM )}\over
{R^*(n;\grM )}}\right|^2&\ll (N^{s/k-1})^{-2}\sum_{n\in {\Bbb N}}\left|
\int_\grm h(\alpha )^se(-n\alpha )d\alpha \right|^2\\
&\ll (N^{s/k-1})^{-2}\int_\grm |h(\alpha )|^{2s}d\alpha .\tag3.1 \endalign $$
When $s\ge {1\over 2}k(\log k+\log \log k+2+o(1))$, the minor arc integral in
(3.1) is $o(N^{2s/k-1})$, and thus it follows that $Z(N)=o(N)$. Thus one may
conclude that almost all integers are sums of $s\sim ({1\over 2}+o(1))k\log k$
positive integral $k$th powers.

The application of Bessel's inequality in (3.1) makes inefficient use of
underlying arithmetic information, and fails, for example, to effectively
estimate the number of values of a polynomial sequence not represented in some
prescribed form. Suppose instead that we define a Fourier series over the
exceptional set itself, namely $K(\alpha )=\sum_ne(n\alpha )$, where the
summation is over $n\in {\Cal Z}(N)$. Since $R^*(n)=0$ for $n\in {\Cal Z}(N)$,
one has $R^*(n;\grm )=-R^*(n;\grM )$, and thus we see that
$$N^{s/k-1}Z(N)\ll \int_\grM h(\alpha )^sK(-\alpha )d\alpha =\left|
\int_\grm h(\alpha )^sK(-\alpha )d\alpha \right| .$$
Applying Schwarz's inequality in combination with Parseval's identity, we
recover the previous consequence of Bessel's inequality via the bound
$$\left| \int_\grm h(\alpha )^sK(-\alpha )d\alpha \right| \le \left(
\int_0^1|K(\alpha )|^2d\alpha \right)^{1/2}\left( \int_\grm |h(\alpha )|^{2s}
d\alpha \right)^{1/2}.\tag3.2$$
However, this formulation permits alternate applications of Schwarz's
inequality or H\"older's inequality. For example, the left hand side of (3.2)
is bounded above by
$$\left( \int_0^1|h(\alpha )^{2t}K(\alpha )^2|d\alpha \right)^{1/2}\left(
\int_\grm |h(\alpha )|^{2s-2t}d\alpha \right)^{1/2},\tag3.3$$
and also by
$$\left( \int_0^1|K(\alpha )|^4d\alpha \right)^{1/4}\left( \int_\grm
|h(\alpha )|^{4s/3}d\alpha \right)^{3/4}.\tag3.4$$
In either case, the diophantine equations underlying the integrals on the
left hand sides of (3.3) and (3.4) contain arithmetic information that can
be effectively exploited whenever the set ${\Cal Z}(N)$ is reasonably thin.

The strategy sketched above has been exploited by Br\"udern, Kawada and Wooley
in a series of papers devoted to additive representation of polynomial
sequences. Typical of the kind of results now available is the conclusion [3]
that almost all values of a given integral cubic polynomial are the sum of six
positive integral cubes. Also, Wooley [30], [31], has derived improved (slim)
exceptional set estimates in Waring's problem when excess variables are
available. For example, write $E(N)$ for the number of integers $n$, with
$1\le n\le N$, for which the anticipated asymptotic formula {\it fails} to hold
 for the number of representations of an integer as the sum of a square
and five cubes of natural numbers. Then in [31] it is shown that $E(N)\ll
N^\epsilon $.

As a final illustration of such ideas, we highlight an application to the
solubility of pairs of diagonal cubic equations. Fix $k=3$, define $h(\alpha )$
 as in \S2, and put $c(n)=\int_0^1|h(\alpha )|^5e(-n\alpha )d\alpha $ for each
$n\in {\Bbb N}$. Br\"udern and Wooley [4] have applied the ideas sketched above
 to estimate the frequency with which large values of $|c(n)|$ occur, and
thereby have shown that, with $\xi $ defined as in the previous section,
$$\sum_{x,y\in {\Cal A}(P,R)}|c(x^3-y^3)|^2=\int_0^1\int_0^1|h(\alpha )^5
h(\beta )^5h(\alpha +\beta )^2|d\alpha d\beta \ll P^{6+\xi +\epsilon }.$$
On noting that $6+\xi <6.25$, cognoscenti will recognise that this twelfth
moment of smooth Weyl sums, in combination with a classical exponential sum
equipped with Weyl's inequality, permits the discussion of pairs of diagonal
cubic equations in $13$ variables via the circle method. The exponent $6+\xi $
improves an exponent $6+2\xi $ previously available for a (different) twelfth
moment. Br\"udern and Wooley [4] establish the following conclusion.

{\bf Theorem 3.1. }\it Suppose that $s\ge 13$, and that $a_i,b_i$ $(1\le i\le s)$ are fixed integers. Then the Hasse
principle holds for the pair of equations
$$a_1x_1^3+\dots +a_sx_s^3=b_1x_1^3+\dots +b_sx_s^3=0.$$
\rm

The condition $s\ge 13$ improves on the previous bound $s\ge 14$ due to
Br\"udern [2], and achieves the theoretical limit of the circle method for this
 problem.

\vskip 8mm

\specialhead \noindent \boldLARGE 4. Arithmetic geometry via descent
\endspecialhead

Let $F({\bold x})\in {\Bbb Z}[x_1,\dots ,x_s]$ be a homogeneous polynomial of degree $d$, and consider the number,
$N(B)$, of integral zeros of the equation $F({\bold x})=0$, with ${\bold x}\in [-B,B]^s$. When $s$ is sufficiently
large in terms of $d$, the circle method shows under modest geometric conditions that $N(B)$ is asymptotic to the
expected product of local
 densities. For fairly general polynomials, the condition on $s$ is as
severe as $s>(d-1)2^d$, though for diagonal equations the methods of \S2 relax
this condition to $s>(1+o(1))d\log d$. However, there is a class of varieties
with small dimension relative to degree, for which the circle method
supplies non-trivial information concerning the density of rational points.
The idea is to apply a descent process in order to interpret points on the
original variety in terms of corresponding points on a new variety, with higher
 dimension relative to degree, more amenable to the circle method.\par

To illustrate this principle, consider a field extension $K$ of ${\Bbb Q}$ of
degree $n$ with associated norm form $N({\bold x})\in {\Bbb Q}[x_1,\dots
,x_n]$. Also, let $l$ and $k$ be natural numbers with $(k,l)=1$, and let
$\alpha $ be a non-zero rational number. Then Heath-Brown and Skorobogatov [7]
descend from the variety $t^l(1-t)^k=\alpha N({\bold x})$ to the associated
variety $aN({\bold u})+bN({\bold v})=z^n$, for suitable integers $a$ and $b$.
The circle method establishes
weak approximation for the latter variety, and thereby it is shown that the
Brauer-Manin obstruction is the only possible obstruction to the Hasse
principle and weak approximation on any smooth projective model of the former
variety. One can artificially construct further examples amenable to the
circle method. For example, if we take linearly independent linear forms
$L_i({\bold x})\in {\Bbb Q} [x_1,\dots ,x_n]$ $(1\le i\le n+r)$, then one can
establish non-trivial lower bounds for the density of rational points on the
variety $z^k=L_1({\bold x})\dots L_{n+r}({\bold x})$ by descending to a variety
 that resembles a system of $r$ diagonal forms of degree $k$, with constrained
varying coefficients. The investigation of such matters will likely provide an
active area of research into the future. In this context we point to work of
Peyre [13], which addresses the interaction between descent and the circle
method in some generality.

\vskip 8mm

\specialhead \noindent \boldLARGE References \endspecialhead \widestnumber \key{999}

\ref
\key1
\by B. J. Birch
\paper \rm Forms in many variables
\jour {\it Proc. Roy. Soc. Ser. A}
\vol \rm 265
\yr 1962
\pages 245--263
\endref

\ref
\key2
\by J. Br\"udern
\paper \rm On pairs of diagonal cubic forms
\jour {\it Proc. London Math. Soc.}
\vol \rm (3) 61
\yr 1990
\pages 273--343
\endref

\ref
\key3
\by J. Br\"udern, K. Kawada and T. D. Wooley
\paper \rm Additive representation in thin sequences, I: Waring's
problem for cubes
\jour {\it Ann. Sci. \'Ecole Norm. Sup.}
\vol \rm (4) 34
\yr 2001
\pages 471--501
\endref

\ref
\key4
\by J. Br\"udern and T. D. Wooley
\paper \rm The Hasse principle for pairs of diagonal cubic equations
\yr to appear
\endref

\ref
\key 5
\by D. R. Heath-Brown
\paper \rm A new form of the circle method, and its application to quadratic
forms
\jour {\it J. Reine Angew. Math.}
\vol \rm 481
\yr 1996
\pages 149--206
\endref

\ref
\key 6
\by D. R. Heath-Brown
\paper \rm Equal sums of three powers
\yr to appear
\endref

\ref
\key 7
\by D. R. Heath-Brown \& A. N. Skorobogatov
\paper \rm Rational solutions of certain equations involving norms
\jour {\it Imperial College preprint}
\yr June 2001
\endref

\ref
\key 8
\by C. Hooley
\paper \rm On nonary cubic forms
\jour {\it J. Reine Angew. Math.}
\vol \rm 386
\yr 1988
\pages 32--98
\endref

\ref
\key 9
\by L.-K. Hua
\paper \rm On Waring's problem
\jour {\it Quart. J. Math. Oxford}
\vol \rm 9
\yr 1938
\pages 199--202
\endref

\ref
\key 10
\by A. A. Karatsuba
\paper \rm Some arithmetical problems with numbers having small prime divisors
\jour {\it Acta Arith.}
\vol \rm 27
\yr 1975
\pages 489--492
\endref

\ref
\key 11
\by Ju. V. Linnik
\paper \rm On the representation of large numbers as sums of seven cubes
\jour \it Mat. Sb.
\vol \rm 12
\yr 1943\pages 218--224
\endref

\ref
\key 12
\by S. T. Parsell
\paper \rm Multiple exponential sums over smooth numbers
\jour \it J. Reine Angew. Math.
\vol \rm 532
\yr 2001
\pages 47--104
\endref

\ref \key 13 \by E. Peyre \paper \rm Torseurs universels et m\'ethode du cercle \inbook {\it Rational points on
algebraic varieties, Progr. Math. 199} \publ Birkh\"auser \yr 2001, 221--274
\endref

\ref
\key 14
\by W. M. Schmidt
\paper \rm The density of integer points on homogeneous varieties
\jour {\it Acta Math.}
\vol \rm 154
\yr 1985
\pages 243--296
\endref

\ref
\key 15
\by R. C. Vaughan
\paper \rm A new iterative method in Waring's problem
\jour {\it Acta Math.}
\vol \rm 162
\yr 1989
\pages 1--71
\endref

\ref
\key 16
\by R. C. Vaughan
\book \it The Hardy-Littlewood Method
\publ Cambridge University Press
\yr 1997
\endref

\ref
\key 17
\by R. C. Vaughan \& T. D. Wooley
\paper \rm Further improvements in Waring's problem, III: eighth powers
\jour \it Philos. Trans. Roy. Soc. London Ser. A
\vol \rm 345
\yr 1993
\pages 385--396
\endref

\ref
\key 18
\by R. C. Vaughan \& T. D. Wooley
\paper \rm Further improvements in Waring's problem, II: sixth powers
\jour \it Duke Math. J.
\vol \rm 76
\yr 1994
\pages 683--710
\endref

\ref
\key 19
\by R. C. Vaughan \& T. D. Wooley
\paper \rm Further improvements in Waring's problem
\jour \it Acta Math.
\vol \rm 174
\yr 1995
\pages 147--240
\endref

\ref
\key 20
\by R. C. Vaughan \& T. D. Wooley
\paper \rm Further improvements in Waring's problem, IV: higher powers
\jour \it Acta Arith.
\vol \rm 94
\yr 2000
\pages 203--285
\endref

\ref
\key 21
\by I. M. Vinogradov
\paper \rm The method of trigonometric sums in the theory of numbers
\jour \it Trav. Inst. Math. Stekloff
\vol \rm 23
\yr 1947
\pages 109
\endref

\ref
\key 22
\by I. M. Vinogradov
\paper \rm On an upper bound for $G(n)$
\jour \it Izv. Akad. Nauk SSSR Ser. Mat.
\vol \rm 23
\yr 1959
\pages 637--642
\endref

\ref
\key 23
\by H. Weyl
\paper \rm \"Uber die Gleichverteilung von Zahlen mod Eins
\jour \it Math. Ann.
\vol \rm 77
\yr 1916
\pages 313--352
\endref

\ref
\key 24
\by T. D. Wooley
\paper \rm Large improvements in Waring's problem
\jour \it Ann. of Math.
\vol \rm (2) 135
\yr 1992
\pages 131--164
\endref

\ref
\key 25
\by T. D. Wooley
\paper \rm On Vinogradov's mean value theorem
\jour \it Mathematika
\vol \rm 39
\yr 1992
\pages 379--399
\endref

\ref
\key 26
\by T. D. Wooley
\paper \rm The application of a new mean value theorem to the fractional
parts of polynomials
\jour \it Acta Arith.
\vol \rm 65
\yr 1993
\pages 163--179
\endref

\ref
\key 27
\by T. D. Wooley
\paper \rm New estimates for smooth Weyl sums
\jour \it J. London Math. Soc.
\vol \rm (2) 51
\yr 1995
\pages 1--13
\endref

\ref
\key 28
\by T. D. Wooley
\paper \rm Breaking classical convexity in Waring's problem: sums of
cubes and quasi-diagonal behaviour
\jour \it Invent. Math.
\vol \rm 122
\yr 1995
\pages 421--451
\endref

\ref
\key 29
\by T. D. Wooley
\paper \rm On exponential sums over smooth numbers
\jour \it J. Reine Angew. Math.
\vol \rm 488
\yr 1997
\pages 79--140
\endref

\ref
\key 30
\by T. D. Wooley
\paper \rm Slim exceptional sets for sums of cubes
\jour \it Canad. J. Math.
\vol \rm 54
\yr 2002
\pages 417--448
\endref

\ref
\key 31
\by T. D. Wooley
\paper \rm Slim exceptional sets in Waring's problem: one
square and five cubes
\jour \it Quart. J. Math.
\vol \rm 53
\yr 2002
\pages 111--118
\endref

\ref
\key 32
\by T. D. Wooley
\paper \rm Sums of three cubes
\jour \it Mathematika
\yr to appear
\endref

\enddocument